\documentclass{amsart}
\usepackage{amsmath}
\usepackage{amsthm}
\usepackage{amsfonts}
\usepackage{amssymb}
\usepackage{graphicx}
\usepackage{fancyhdr}%
\usepackage{xcolor}%
\usepackage[titletoc,title]{appendix}
\usepackage{combelow}
\usepackage{url}

\title{A Note on Fractional Sobolev inequalities revisited}
\author{Giovanni Leoni}
\address{Department of Mathematical Sciences, Carnegie Mellon University, Pittsburgh PA 15213-3890, USA}
\email{giovanni@andrew.cmu.edu}
\author{Asher Wheaton}
\address{Department of Mathematical Sciences, Carnegie Mellon University, Pittsburgh PA 15213-3890, USA}
\email{awheaton@andrew.cmu.edu}
\date{}

\newtheorem{theorem}{Theorem}[section]

\newtheorem{definition}{Definition}[section]

\begin{document}

\maketitle

\begin{abstract}
    In this paper, we extend the analysis performed in \cite{Nguyen} by Nguyen et al. of homogeneous Sobolev inequalities. We combine their proof with the first author's adaptation in \cite{Leoni} to provide a short proof through interpolation for the critical case.
\end{abstract}

\section{Introduction}
In \cite{Nguyen}, Nguyen et al. used a novel application of ordinary and sharp Hardy-Littlewood maximal functions to prove that for $s_1-\frac{N}{p_1}=s_2-\frac{N}{p_2}$, $s_2<s_1$ the embedding
\begin{align*}
    \dot{W}^{s_1,p_1}(\mathbb{R}^N)\hookrightarrow \dot{W}^{s_2,p_2}(\mathbb{R}^N)
\end{align*} 
of homogeneous Sobolev spaces holds. This result was derived by applying the inequality
\begin{flalign*}
    \int_{\mathbb{R}^N}\int_{\mathbb{R}^N} \frac{|u(x+h)-u(x)|^q}{\|h\|^{N+\sigma q}}dhdx \lesssim&&
\end{flalign*}
\begin{flalign*}
    &&\int_{\mathbb{R}^N} (M(u)(x))^{(1-\theta)q}\left(\sup_{r>0}\frac{1}{r^N}\int_{B(0,r)}|u(x+z)-u(x)| \frac{dz}{\|z\|^{\sigma}} \right)^{\theta q}dx,
\end{flalign*}
where $1\leq q<\infty$ and $0<\theta,\sigma<1$, in combination with the Sobolev-Gagliardo-Nirenberg embedding. This method relies on the intrinsic Gagliardo semi-norm, 
$$|u|_{W^{s,p}(\mathbb{R}^N)}=\left(\int_{\mathbb{R}^N}\int_{\mathbb{R}^N}\frac{|u(x+h)-u(x)|^p}{\|h\|^{N+sp}}dhdx\right)^{1\over p},$$
where $1\leq p, q < \infty$ and $0<s<1.$ The method offers a simpler alternative to classical proofs, including those from \cite{Triebel} and \cite{Solonnikov}. The former is presented as Theorem 17.85 of \cite{LeoniNonFrac} and uses the Littlewood-Paley decomposition to prove the critical and supercritical cases, that is, when $s_1-\frac{N}{p_1}= 0$ and $s_1-\frac{N}{p_1}>0$ respectively. The latter is presented as Theorem 7.35 in \cite{Leoni} and uses anisotropic $L^p$ spaces to prove the general case in Besov spaces.

Although the proof from Nguyen et al. is simpler, their analysis cannot cover the critical and supercritical cases when $s_1<1$ and $s_2<1$. In these settings, standard embeddings are not sufficiently compatible with maximal functions for an identical approach. 

The first author extended the result of Nguyen et al. in Theorem 7.34 of \cite{Leoni} by modifying the Hardy-Littlewood maximal function. He proved the inequality 
\begin{flalign*}
    \int_{\mathbb{R}^N}\left(\frac{1}{\|h\|^N}\int_{B(0,2\|h\|)}|u(x+z)-u(x)|dz\right)^{p_2}\frac{dh}{\|h\|^{N+s_2p_2}}\lesssim&&
\end{flalign*}
\begin{flalign*}
    &&|u|_{C^{0,\alpha}(\mathbb{R}^N)}^{p_2-p_1}\left(\sup_{r>0}\frac{1}{r^N}\int_{B(0,r)}|u(x+z)-u(x)|\frac{dz}{\|z\|^{s_1}}\right)^{p_1},
\end{flalign*}
where $0<s_2<s_1<1$, $1\leq p_1,p_2 \leq \infty$, $s_1-\frac{N}{p_1}=s_2-\frac{N}{p_2}=\alpha>0$. After an application of Morrey's embedding theorem, this handles the supercritical case. Nevertheless, the first author resorted to Solonnikov's long and classical method from \cite{Solonnikov} to prove the critical case in Theorem 7.35. An alternative proof akin to Nguyen et al.'s result may be available through the space $\operatorname{BMO}(\mathbb{R}^N)$, although no such proof is forthcoming. Instead, in this paper, we succinctly prove the critical case using standard results of interpolation and the short theorems taken from both \cite{Leoni} and \cite{Nguyen}. Our result is as follows.

\begin{theorem} \label{main}
    For $0<s_2<s_1<1$, $1\leq p_1,p_2< \infty$ where $s_1-\frac{N}{p_1}=s_2-\frac{N}{p_2}=0$,
    \begin{align}
        \dot{W}^{s_1,p_1}(\mathbb{R}^N)\hookrightarrow \dot{W}^{s_2,p_2}(\mathbb{R}^N).
    \end{align}
\end{theorem}

\section{Preliminaries}

In this paper, we write $A(x)\lesssim B(x)$ if for some constant $c>0$ independent of x, $A(x)\leq cB(x)$ for all $x\in dom(A)$.

Combining the aforementioned results of \cite{Leoni} and \cite{Nguyen}, the statement is proven in the noncritical case.
\begin{theorem} \label{noncritical}
    For $0< s_2<s_1<1$, $1\leq p_1,p_2 < \infty$ where $s_1-\frac{N}{p_1}=s_2-\frac{N}{p_2}\neq0$,
    \begin{align}
        \dot{W}^{s_1,p_1}(\mathbb{R}^N)\hookrightarrow \dot{W}^{s_2,p_2}(\mathbb{R}^N).
    \end{align}
\end{theorem}
We make use of several classical results on interpolation, which are extended to interpolation on semi-normed spaces. A more comprehensive discussion of the subject can be found in \cite{LeoniNonFrac} and \cite{Stevenson}, although we present a brief overview of the relevant terms and theorems in this paper. 


\begin{definition}
Two semi-normed spaces $(X_0, |\cdot|_{X_0})$ and $(X_1, |\cdot|_{X_1})$ are compatible if both embed into the same larger vector space $X$.
\end{definition}

Both the intersection and the sum of compatible spaces inherit natural semi-norm structures. The real interpolation spaces defined below are particular examples of intermediate spaces between $X_0\cap X_1$ and $X_0 + X_1.$

\begin{definition}
For compatible semi-normed spaces $(X_0, |\cdot|_{X_0})$ and $(X_1, |\cdot|_{X_1})$ and any $x\in X_0+X_1$, define the semi-norm $K(x,t;X_0,X_1)=\inf(|x_0|_{X_0}+t|x_1|_{X_1})$, where the infimum is taken over all $x_0 \in X_0,x_1 \in X_1$ with $x_0+x_1=x.$ Where $X_0$ and $X_1$ are apparent, we write $K(x,t)$ instead.
\end{definition}
\begin{definition}[K Method of Interpolation] \label{KMethod}
For compatible semi-normed spaces $(X_0, |\cdot|_{X_0})$ and $(X_1, |\cdot|_{X_1})$, $0<\theta<1$, and $1\leq p \leq \infty$, define 
\begin{align} \label{star}
    |x|_{(X_0,X_1)_{\theta,p}}=\left(\int_0^\infty (K(x,t))^p\frac{dt}{t^{1+\theta p}}\right)^{\frac{1}{p}}
\end{align} 
if $p<\infty$, and
\begin{align*}
    |x|_{(X_0,X_1)_{\theta,p}}=\sup_{t>0}t^{-\theta}K(x,t)
\end{align*}
otherwise. $(X_0,X_1)_{\theta,p}$ is the set of all $x\in X_0+X_1$ such that this value is finite.
\end{definition}

Spaces generated by the K-method have several useful properties, the most important of which is the following.
\begin{theorem}[Reiteration] \label{reit}
    Let $(X_0, |\cdot|_{X_0})$ and $(X_1, |\cdot|_{X_1})$ be compatible complete semi-normed spaces with complete intersection. For $0\leq s_0< s_1\leq 1$, $0<\theta<1$, $s=\theta s_1 + (1-\theta)s_0$, $1\leq q_0, q_1, q_2\leq\infty$,
    \begin{align}
        (X_{s_0,q_0},X_{s_1,q_1})_{\theta,q_2}=(X_0,X_1)_{s,q_2}
    \end{align}
    with equivalence of seminorms, where 
    \begin{align}
        X_{s_i,q_i}=\begin{cases} X_{s_i} & s_i\in \{0,1\}, \\ (X_0,X_1)_{s_i,q_i} & s_i \not\in \{0,1\}.\end{cases}
    \end{align}
\end{theorem}
The reiteration theorem is proved in \cite{Stevenson}. Although it requires a few tricks, the method is largely long algebraic manipulation. We also use a few standard inclusions of interpolation spaces.
\begin{theorem} \label{embedup}
    Let $(X_0, |\cdot|_{X_0}), (X_1, |\cdot|_{X_1})$ be an admissible pair of semi-normed spaces, $0<\theta<1$, $1\leq p\leq q < \infty$. In this setting,
    \begin{align}
        (X_0,X_1)_{\theta,p} \hookrightarrow (X_0,X_1)_{\theta,q} \hookrightarrow (X_0,X_1)_{\theta,\infty}.
    \end{align}
\end{theorem}
\begin{proof}
We present a proof taken from \cite{LeoniNonFrac}. As $K(x,\cdot)$ is increasing, for $t>0,$ by (\ref{star}),
\begin{align*}
    t^{-\theta}K(x,t)&\lesssim\left(\int_t^\infty (K(x,t))^q\frac{dz}{z^{1+\theta q}}\right)^{1/q} \\
    &\leq\left(\int_t^\infty (K(x,z))^q\frac{dz}{z^{1+\theta q}}\right)^{1/q} \\
    &\leq|x|_{(X_0,X_1)_{\theta,q}}.
\end{align*}
Taking the supremum of the left side yields the second embedding. 
Further, 
\begin{align*}
    |x|_{(X_0, X_1)_{\theta, q}}^q&=\int_0^\infty (K(x,t))^q\frac{dt}{t^{1+\theta q}}\\
    &=\int_0^\infty (t^{-\theta}K(x,t))^{q-p}(K(x,t))^p\frac{dt}{t^{1+\theta p}}\\
    &\lesssim|x|_{(X_0, X_1)_{\theta, \infty}}^{q-p}\int_0^\infty (K(x,t))^p\frac{dt}{t^{1+\theta p}}\\
    &\lesssim|x|_{(X_0, X_1)_{\theta, p}}^q,
\end{align*}
where the final inequality uses the second embedding applied to $(X_0, X_1)_{\theta, p}$.
\end{proof}

\begin{theorem} \label{bdd}
    Let $(X_0, |\cdot|_{X_0}),(X_1, |\cdot|_{X_1})$ and $(Y_0, |\cdot|_{Y_0}),(Y_1, |\cdot|_{Y_1})$ be two pairs of compatible semi-normed spaces, and let $T:X_0+X_1\rightarrow Y_0+Y_1$ be linear. Suppose $T:X_i\rightarrow Y_i$ and $|T(x_i)|_{Y_i}\lesssim|x_i|_{X_i}$ for all $x_i\in X_i$, $i=0,1$. For all $0<\theta<1$, $1\leq p \leq \infty$, $|T(x)|_{(Y_0,Y_1)_{\theta,p}}\lesssim|x|_{(X_0,X_1)_{\theta,p}}$ for all $x \in (X_0, X_1)_{\theta,p}$.
\end{theorem}
\begin{proof} We present an adapted proof from \cite{LeoniNonFrac}. Let $x_0\in X_0$, $x_1\in X_1$ be such that $x_0+x_1=x$. By linearity,
\begin{align*}
    K(T(x),t;Y_0, Y_1)\leq |T(x_0)|_{Y_0}+t|T(x_1)|_{Y_1}\lesssim |x_0|_{X_0}+t|x_1|_{X_1}
\end{align*}
Taking the infimum yields $K(T(x),t; Y_0,Y_1)\lesssim K(x,t;X_0,X_1)$. The result follows.
\end{proof}

We conclude our preliminaries with a specific application of the K-method to Sobolev spaces.

\begin{theorem} \label{interp}
    For $1 \leq p < \infty$, $0<\theta<1$, 
    \begin{align}
        (L^p\left(\mathbb{R}^N\right), \dot{W}^{1,p}\left(\mathbb{R}^N\right))_{\theta,p} = \dot{W}^{\theta,p}\left(\mathbb{R}^N\right).
    \end{align}
\end{theorem}

This result is proven in \cite{Leoni} using mollification to decompose an arbitrary function. Note that $\dot{W}^{1,p}(\mathbb{R}^N)\cap L^p(\mathbb{R}^N)=W^{1,p}(\mathbb{R}^N)$ is a Banach space. Thus, we may apply the reiteration theorem to Sobolev spaces.

\section{Proof}

\begin{proof}[Proof of Theorem \ref{main}]
    Observe that $1<p_1=\frac{N}{s_1}<p_2=\frac{N}{s_2}<\infty.$ Let $0<d_1<s_1<t_1<1$ and let $0<\theta<1$ such that $\theta t_1+(1-\theta)d_1=s_1$. Define $d_2=d_1+\frac{N}{p_2}-\frac{N}{p_1}$ and $t_2=t_1+\frac{N}{p_2}-\frac{N}{p_1}$. Because \begin{equation*}
        d_1-\frac{N}{p_1}=d_2-\frac{N}{p_2}=d_1-s_1<0, \\
        t_1-\frac{N}{p_1}=t_2-\frac{N}{p_2}=t_1-s_1>0,
    \end{equation*} we may apply Theorems \ref{reit} and \ref{interp} to write 
    \begin{align*}
    \dot{W}^{s_1,p_1}(\mathbb{R}^N)=(\dot{W}^{d_1,p_1}(\mathbb{R}^N),\dot{W}^{t_1,p_1}(\mathbb{R}^N))_{\theta,p_1}.
    \end{align*}
    By taking $d_1, t_1$ close to $s_1$, we ensure $0<d_2,t_2<1$, since $0<s_2=s_1+\frac{N}{p_2}-\frac{N}{p_1}<1$. By Theorems \ref{noncritical} and \ref{bdd}, 
    \begin{align*}
        (\dot{W}^{d_1,p_1}(\mathbb{R}^N),\dot{W}^{t_1,p_1}(\mathbb{R}^N))_{\theta,p_1} \hookrightarrow (\dot{W}^{d_2,p_2}(\mathbb{R}^N),\dot{W}^{t_2,p_2}(\mathbb{R}^N))_{\theta,p_1},
    \end{align*}
    where we use the identity map as our linear transformation. By Theorem \ref{embedup} 
    \begin{align*}
        (\dot{W}^{d_2,p_2}(\mathbb{R}^N),\dot{W}^{t_2,p_2}(\mathbb{R}^N))_{\theta,p_1} \hookrightarrow (\dot{W}^{d_2,p_2}(\mathbb{R}^N),\dot{W}^{t_2,p_2}(\mathbb{R}^N))_{\theta,p_2}.
    \end{align*}
    A final application of Theorems \ref{reit} and \ref{interp} completes the proof, since
    \begin{align*}
        s_2=s_1+\frac{N}{p_2}-\frac{N}{p_1}=\theta t_2 +(1-\theta)d_2.
    \end{align*}
\end{proof}

\section*{Acknowledgments}
Asher Wheaton was supported by the summer research support of NSF grant  DMS \#2508400
This paper was written without using AI.
\bibliographystyle{abbrv}
\bibliography{sources}

\end{document}